\documentclass[a4paper,12pt]{article}

\usepackage[all]{xy}
\usepackage[T1]{fontenc}
\usepackage[dvips]{graphicx}
\usepackage{amsmath,amssymb,amsfonts,amsthm,latexsym,mathrsfs,textcomp,verbatim, enumerate, soul, color}
\xyoption{web}
\begin{document}

\title{A general method for building reflections}

\author{Olivia Caramello \vspace{3 mm}\\ {\small DPMMS, University of Cambridge,}\\{\small Wilberforce Road, Cambridge CB3 0WB, U.K.}\\{\small O.Caramello@dpmms.cam.ac.uk}\thanks{The author gratefully acknowledges the support of a Research Fellowship from Jesus College, Cambridge (U.K.)}}

\date{December 15, 2011}

\maketitle

\def\Monthnameof#1{\ifcase#1\or
   January\or February\or March\or April\or May\or June\or
   July\or August\or September\or October\or November\or December\fi}
\def\today{\number\day~\Monthnameof\month~\number\year}

%
%
%
\def\pushright#1{{
   \parfillskip=0pt            
   \widowpenalty=10000         
   \displaywidowpenalty=10000  
   \finalhyphendemerits=0      
  %
   \leavevmode                 
   \unskip                     
   \nobreak                    
   \hfil                       
   \penalty50                  
   \hskip.2em                  
   \null                       
   \hfill                      
   {#1}                        
  %
   \par}}                      

\def\qed{\pushright{$\square$}\penalty-700 \smallskip}

\newtheorem{theorem}{Theorem}[section]

\newtheorem{proposition}[theorem]{Proposition}

\newtheorem{scholium}[theorem]{Scholium}

\newtheorem{lemma}[theorem]{Lemma}

\newtheorem{corollary}[theorem]{Corollary}

\newtheorem{conjecture}[theorem]{Conjecture}

\newenvironment{proofs}%
 {\begin{trivlist}\item[]{\bf Proof }}%
 {\qed\end{trivlist}}

  \newtheorem{rmk}[theorem]{Remark}
\newenvironment{remark}{\begin{rmk}\em}{\end{rmk}}

  \newtheorem{rmks}[theorem]{Remarks}
\newenvironment{remarks}{\begin{rmks}\em}{\end{rmks}}

  \newtheorem{defn}[theorem]{Definition}
\newenvironment{definition}{\begin{defn}\em}{\end{defn}}

  \newtheorem{eg}[theorem]{Example}
\newenvironment{example}{\begin{eg}\em}{\end{eg}}

  \newtheorem{egs}[theorem]{Examples}
\newenvironment{examples}{\begin{egs}\em}{\end{egs}}


\mathcode`\<="4268  
\mathcode`\>="5269  
\mathcode`\.="313A  
\mathchardef\semicolon="603B 
\mathchardef\gt="313E
\mathchardef\lt="313C

\newcommand{\imp}
 {\!\Rightarrow\!}

\newcommand{\biimp}
 {\!\Leftrightarrow\!}

\newcommand{\bjg}
 {\mathrel{{\dashv}\,{\vdash}}}

\newcommand{\cod}
 {{\rm cod}}

\newcommand{\dom}
 {{\rm dom}}

\newcommand{\epi}
 {\twoheadrightarrow}

\newcommand{\Ind}[1]
 {{\rm Ind}\hy #1}

\newcommand{\mono}
 {\rightarrowtail}

\newcommand{\nml}
 {\triangleleft}

\newcommand{\ob}
 {{\rm ob}}

\newcommand{\op}
 {^{\rm op}}

\newcommand{\pepi}
 {\rightharpoondown\kern-0.9em\rightharpoondown}

\newcommand{\pmap}
 {\rightharpoondown}

\newcommand{\Set}
 {{\bf Set }}

\newcommand{\Sh}
 {{\bf Sh}}

\newcommand{\sh}
 {{\bf sh}}

\newcommand{\Sub}
 {{\rm Sub}}

\begin{abstract}
We establish a general method for generating reflections between categories. We then apply our technique to generate adjunctions starting from geometric morphisms between Grothendieck toposes; as particular cases, we recover various well-known Stone-type adjunctions and establish several new ones.
\end{abstract}

\section{Introduction}

Adjunction is a fundamental relationship between pairs of categories. An illuminating point of view on this notion is provided by the concept of comma category, originally introduced by F. W. Lawvere in his Ph.D. thesis \cite{lawvere}; indeed, the fact that two functors $F:{\cal A}\to {\cal B}$ and $G:{\cal B}\to {\cal E}$ between a given pair of categories are adjoint to each other can be expressed as the existence of an equivalence between the two comma categories $(F\downarrow 1_{{\cal B}})$ and $(1_{\cal A} \downarrow G)$, and an adjunction $F\dashv G$ can be obtained by composing a coreflection between $\cal A$ and $(F\downarrow 1_{{\cal B}})$ with the equivalence $(F\downarrow 1_{{\cal B}})\simeq (1_{\cal A}\downarrow G)$ and then with a reflection between $(1_{\cal A}\downarrow G)$ and $\cal B$. The comma category thus acts as a `bridge' object which condenses the information about the adjunction (via its two different representations) and connects the two categories with each other. 

The aspect of the comma category construction which constitutes the inspiration for the present paper is the fact that this construction shows that the relationship between two categories $\cal A$ and $\cal B$ related by a pair of adjoint functors $F$ and $G$ may well be best understood from the point of view of a third category, namely the comma category $(F\downarrow 1_{{\cal B}})\simeq (1_{\cal A}\downarrow G)$, into which the two categories $\cal A$ and $\cal B$ (canonically) embed. 

This naturally leads to the idea of a general method for building adjunctions between a given pair of categories starting from a pair of functors from each of the two categories into a third one, together with some relationships between them. In this paper, we show that this intuition can be materialized in a precise technical sense, by providing in section \ref{general} a general method for building reflections between categories starting from data of that kind. As shown in section \ref{completeness}, our method is \emph{complete}, in the sense that \emph{any} reflection between categories can be obtained as an application of it. 

The interest of our method lies in its inherent technical flexibility; indeed, it happens very often in practice that two different categories are best understood in relationship with each other from the point of view of a third category to which both are related (cf. for example \cite{OC10} for an explanation of the sense in which Grothendieck toposes can act as `bridges' connecting different mathematical theories with each other), and our method allows us to establish reflections between a given pair of categories starting from relations between the `realizations' of the two categories at the level of the `bridge category' to which both categories map.

In \cite{OC11} we showed that many `Stone-type' dualities or equivalences between categories of preorders and categories of posets, locales and topological spaces can be naturally interpreted as arising from the process of appropriately `functorializing' categorical equivalences between toposes; in that paper, we also established various adjunctions between categories of these kinds which extend the given dualities or equivalences. Amongst other things, in section \ref{geomrefl} of this paper, we show as applications of our general method that all the `Stone-type' reflections obtained in that paper can be seen as applications of our method for building reflection in the context where the `bridge category' is some category of toposes or of toposes paired with points (as defined in \cite{OC11}).      

\newpage
 
\section{The general method}\label{general}

In this section we describe our general method for generating reflections between categories. 

The set of data we shall work with consists of two categories $\cal H$ and $\cal K$, a category $\cal U$, two functors $I:{\cal H}\to {\cal U}$, $J:{\cal K}\to {\cal U}$ and two binary relations $R$ and $S$ on $Ob({\cal H})\times Ob({\cal K})$. In addition to this, we suppose to have, for every $({\cal C}, {\cal D})\in R$, an arrow $\xi_{({\cal C}, {\cal D})}:I({\cal C})  \to J({\cal D})$ in $\cal U$ and for every $({\cal C}, {\cal D})\in S$ an arrow $\chi_{({\cal C}, {\cal D})}:J({\cal D}) \to I({\cal C})$ in $\cal U$, and two functions $Z:R\to S$ and $W:S\to R$ such that $Z$ keeps the second component fixed and $W$ keeps the first component fixed. 

Let us denote by $\pi^{{\cal H}}_{R}$ and $\pi^{{\cal K}}_{R}$ respectively the canonical projections $R\hookrightarrow {\cal H}\times {\cal K}\to {\cal H}$ and $R\hookrightarrow {\cal H}\times {\cal K}\to {\cal K}$; similarly, we denote by $\pi^{{\cal H}}_{S}$ and $\pi^{{\cal K}}_{S}$ respectively the canonical projections $S\hookrightarrow {\cal H}\times {\cal K}\to {\cal H}$ and $S\hookrightarrow {\cal H}\times {\cal K}\to {\cal K}$. We impose the following hypotheses.

\begin{enumerate}
\item For any $({\cal C}, {\cal D})\in R$, $Z(({\cal C}, {\cal D}))=(\pi^{{\cal H}}_{S}(Z(({\cal C}, {\cal D}))), {\cal D})$; we require the composite
\[
\chi_{Z({\cal C}, {\cal D})}\circ \xi_{({\cal C}, {\cal D})}: I({\cal C}) \to J(\pi^{{\cal H}}_{S}(Z(({\cal C}, {\cal D}))))
\]
to be induced by a (canonically chosen) isomorphism in $\cal H$ 
\[
{{\epsilon'_{({\cal C}, {\cal D})}}}^{-1}: {\cal C} \to \pi^{{\cal H}}_{S}(Z(({\cal C}, {\cal D}))),
\]
as in the following diagram:
\[  
\xymatrix {
I({\cal C}) \ar[r]^{\xi_{({\cal C}, {\cal D})}} \ar[dr]^{I({{\epsilon'_{({\cal C}, {\cal D})}}}^{-1})} & J({\cal D}) \ar[d]^{\chi_{Z({\cal C}, {\cal D})}}\\
& I(\pi^{{\cal H}}_{S}(Z(({\cal C}, {\cal D}))))}
\]

\item For any $({\cal C}, {\cal D})\in S$, $W(({\cal C}, {\cal D}))=({\cal C}, \pi^{{\cal K}}_{R}(W(({\cal C}, {\cal D}))))$; we require the composite
\[
\xi_{W({\cal C}, {\cal D})}\circ \chi_{({\cal C}, {\cal D})}: J({\cal D})\to J(\pi^{{\cal K}}_{R}(W(({\cal C}, {\cal D}))))
\]
to be induced by a (canonically chosen) morphism
\[
{\eta'_{({\cal C}, {\cal D})}}:{\cal D} \to \pi^{{\cal K}}_{R}(W(({\cal C}, {\cal D})))
\]
in $\cal K$, as in the following diagram:
\[  
\xymatrix {
J({\cal D}, K_{{\cal D}}) \ar[r]^{\chi_{({\cal C}, {\cal D})}} \ar[dr]^{I({\eta'_{({\cal C}, {\cal D})}})} & I({\cal C}) \ar[d]^{\xi_{W({\cal C}, {\cal D})}}\\
& J(\pi^{{\cal K}}_{R}(W(({\cal C}, {\cal D}))))}
\]

\item for any $({\cal C}, {\cal D})\in R$, the arrow
\[
\xi_{({\cal C}, {\cal D})}:I({\cal C})\to J({\cal D})
\]
is an isomorphism in $\cal U$. [Note that, since $\epsilon'_{({\cal C}, {\cal D})}$ is an isomorphism, this implies that 
\[
\chi_{Z({\cal C}, {\cal D})}:J({\cal D}) \to I(\pi^{{\cal H}}_{S}(Z(({\cal C}, {\cal D}))))
\]
is an isomorphism as well.]

\item For any $({\cal C}, {\cal D})\in R$, $\eta'_{Z(({\cal C}, {\cal D}))}$ is an isomorphism.
\end{enumerate} 
 
Let us define two categories $\tilde{R}$ and $\tilde{S}$, as follows. 

The objects of $\tilde{R}$ are the elements of $R$ while the arrows $({\cal C}, {\cal D})\to ({\cal C}', {\cal D}')$ are the pairs $(u,v)$ where 
\[
u:{\cal C}\to {\cal C}',
\]
\[
v:{\cal D}\to {\cal D}'
\]
are arrows respectively in the categories $\cal H$, $\cal K$ and such that the following square commutes:

\[  
\xymatrix {
I({\cal C}') \ar[d]^{\xi_{({\cal C}', {\cal D}')}} \ar[r]^{I(u)} & I({\cal C})  \ar[d]^{\xi_{({\cal C}, {\cal D})}}\\
J({\cal D}') \ar[r]^{J(v)} & J({\cal D}) }
\]   

We will occasionally write $(u,v,z)$ for $(u,v)$, where $z$ is the arrow 
\[
\pi^{{\cal H}}_{S}(Z(({\cal C}, {\cal D}))) \to \pi^{{\cal H}}_{S}(Z(({\cal C}', {\cal D}')))
\]
in $\cal H$ given by the factorization of $u$ across the isomorphisms $\epsilon'_{({\cal C}, {\cal D})}$ and $\epsilon'_{({\cal C}', {\cal D}')}$. 

The composition of arrows in $\tilde{R}$ is defined as the composition of the functors forming the various components.

Similarly, we define the category $\tilde{S}$. The objects of $\tilde{S}$ are the elements of $S$ while the arrows $({\cal C}, {\cal D})\to ({\cal C}', {\cal D}')$ are the triples $(z,v,w)$, where 
\[
v:{\cal D}\to {\cal D}',
\]
\[
z:{\cal C} \to {\cal C}'
\]
and 
\[
w:\pi^{{\cal K}}_{R}(W(({\cal C}, {\cal D}))) \to \pi^{{\cal K}}_{R}(W(({\cal C}', {\cal D}')))
\]
are morphisms respectively in the categories $\cal K$, $\cal H$ and $\cal K$ such that the two squares in the following diagram commute:

\[  
\xymatrix {
J({\cal D}') \ar[d]^{\chi_{({\cal C}', {\cal D}')}} \ar[r]^{J(v)} & J({\cal D})  \ar[d]^{\chi_{({\cal C}, {\cal D})}}\\
I({\cal C}') \ar[d]^{\xi_{({\cal C}', {\cal D}')}} \ar[r]^{I(z)} & I({\cal C})  \ar[d]^{\xi_{({\cal C}, {\cal D})}}\\
J(\pi^{{\cal K}}_{R}(W(({\cal C}', {\cal D}')))) \ar[r]^{J(w)} & J(\pi^{{\cal K}}_{R}(W(({\cal C}, {\cal D}))))}.
\]   

Composition of arrows in $\tilde{S}$ is defined componentwise as the composition of the functors forming the various components.

Let us now define functors $\tilde{Z}:\tilde{R}\to \tilde{S}$ and $\tilde{W}:\tilde{S} \to \tilde{R}$, which extend the functions $Z:R\to S$ and $W:S\to R$. 

For any object $({\cal C}, {\cal D})$ of $\tilde{R}$, we set $\tilde{Z}(({\cal C}, {\cal D}))=Z(({\cal C}, {\cal D}))$ and for any arrow $(u,v,z):({\cal C}, {\cal D}) \to ({\cal C}', {\cal D}')$ in $\tilde{R}$ we set $\tilde{Z}((u,v,z))$ equal to the triple $(z, v, w)$, where $w:\pi^{{\cal K}}_{R}(W(Z({\cal C}, {\cal D})) \to \pi^{{\cal K}}_{R}(W(Z({\cal C}', {\cal D})')$ is the only arrow in $\cal K$ making the following diagram commute (recall that by our hypotheses $\eta'_{Z({\cal C}, {\cal D})}$ and $\eta'_{Z({\cal C}', {\cal D}')}$ are isomorphisms). 

\[  
\xymatrix {
{\cal D} \ar[d]^{\eta'_{Z({\cal C}, {\cal D})}} \ar[r]^{v} & {\cal D}' \ar[d]^{\eta'_{Z({\cal C}', {\cal D}')}}\\
\pi^{{\cal K}}_{R}(W(Z({\cal C}, {\cal D}))) \ar[r]^{w}   & \pi^{{\cal K}}_{R}(W(Z({\cal C}', {\cal D})'))  }
\]   

To show that $(z,v,w)$ is an arrow in $\tilde{S}$ we have to check that the diagram

\[  
\xymatrix {
J({\cal D}') \ar[d]^{\chi_{\cal D}'} \ar[r]^{J(v)} & J({\cal D})  \ar[d]^{\chi_{({\cal C}, {\cal D})}}\\
I(\pi^{{\cal H}}_{S}(Z(({\cal C}', {\cal D}')))') \ar[d]^{\xi_{(\pi^{{\cal H}}_{S}(Z(({\cal C}, {\cal D}))))}} \ar[r]^{I(z)} & I({(\pi^{{\cal H}}_{S}(Z(({\cal C}, {\cal D}))))})  \ar[d]^{\xi_{{(\pi^{{\cal H}}_{S}(Z(({\cal C}, {\cal D}))))}}}\\
J(\pi^{{\cal K}}_{R}(W(Z(({\cal C}', {\cal D}'))))) \ar[r]^{J(w)} & J(\pi^{{\cal K}}_{R}(W(Z(({\cal C}, {\cal D})))))}.
\]  

commutes. 

The top square commutes since $(u,v,z)$ is an arrow $({\cal C}, {\cal D})\to ({\cal C}', {\cal D}')$ in $\tilde{R}$. It remains to prove the commutativity of the bottom square; but this is equivalent to the commutativity of the outer rectangle, since $\chi_{\cal D}'$ and $\chi_{\cal D}$ are isomorphisms (cf. our assumptions above), and to show this it suffices to invoke the commutativity of the square defining $w$.

This completes the proof that our assignment defines a functor $\tilde{Z}:\tilde{R}\to \tilde{S}$. Let us now turn to the definition of the functor $\tilde{W}:\tilde{S}\to \tilde{R}$. For any $({\cal C}, {\cal D})$ in $\tilde{S}$, we set $\tilde{W}(({\cal C}, {\cal D}))$ equal to $W(({\cal C}, {\cal D}))$ and for any arrow $(z,v,w):({\cal C}, {\cal D})\to ({\cal C}', {\cal D}')$ in $\tilde{S}$, we set $\tilde{W}((z,v,w))=(z,w)$. This is clearly an arrow in $\tilde{R}$ and hence this assignment actually defines a functor, as required.

We are now ready to state our main result.

\begin{theorem}
Under the hypotheses specified above, the functors $\tilde{Z}:\tilde{R}\to \tilde{S}$ and $\tilde{W}:\tilde{S}\to \tilde{R}$ are adjoint to each other; in fact, they yield a reflection in which $\tilde{Z}$ is the right adjoint and $\tilde{W}$ is the left adjoint. 
\end{theorem}

\begin{proofs}
We shall establish this adjunction by giving its unit and counit and checking that they are natural transformations satisfying the triangular identities. We take as counit the transformation $\epsilon$ sending each object $({\cal C}, {\cal D})$ of $\tilde{R}$ the pair $\epsilon_{({\cal C}, {\cal D})}:=(\epsilon'_{({\cal C}, {\cal D})}, {\eta'_{Z({\cal C}, {\cal D})}}^{-1})$, regarded as an arrow $\tilde{W}(\tilde{Z}(({\cal C}, {\cal D})))=(\pi^{{\cal H}}_{S}(Z(({\cal C}, {\cal D}))), \pi^{{\cal K}}_{R}(W(Z({\cal C}, {\cal D}))))\to ({\cal C}, {\cal D})$ in $\tilde{R}$. This is indeed an arrow in $\tilde{R}$ by the commutativity of the following diagram.

\[  
\xymatrix {
I({\cal C}) \ar[d]^{\xi_{({\cal C}, {\cal D})}}  & & & & \ar[llll]^{I({\epsilon'_{({\cal C}, {\cal D})}})} I(\pi^{{\cal H}}_{S}(Z(({\cal C}, {\cal D})))) \ar[d]^{\xi_{W(Z(({\cal C}, {\cal D})))}} \\
J({\cal D}) \ar[urrrr]_{\chi_{Z({\cal C}, {\cal D})}}  & & & & J(\pi^{{\cal K}}_{R}(W(Z({\cal C}, {\cal D})))) \ar[llll]^{J({\eta'_{Z({\cal C}, {\cal D})}}^{-1})}}
\]   

Note that $\epsilon_{({\cal C}, {\cal D})}=({\epsilon'_{Z({\cal C}, {\cal D})}}, {\eta'_{({\cal C}, {\cal D})}}^{-1}, z)$ where $z$ is the unique morphism making the following diagram commute:
 
\[  
\xymatrix {
{\cal C} & & \pi_{S}^{\cal H}(Z({\cal C}, {\cal D})) \ar[ll]^{{\epsilon'_{({\cal C}, {\cal D})}}} \\
\pi_{S}^{\cal H}(Z({\cal C}, {\cal D})) \ar[u]^{{\epsilon'_{({\cal C}, {\cal D})}}}  & & \ar[ll]^{z} \ar[u]^{{\epsilon'_{W(Z({\cal C}, {\cal D}))}}} \pi_{S}^{\cal H}(Z(W(Z(({\cal C}, {\cal D})))))}
\]
that is, $z={\epsilon'_{W(Z({\cal C}, {\cal D}))}}$. So we have 
\[
\epsilon_{({\cal C}, {\cal D})}=({\epsilon'_{({\cal C}, {\cal D})}}, {\eta'_{Z({\cal C}, {\cal D})}}^{-1}, {\epsilon'_{W(Z({\cal C}, {\cal D}))}}).
\] 

The transformation $\epsilon$ is natural because each of its components is arise from compositions of (morphisms which are naturally isomorphic to) the morphisms $\xi$ and $\chi$ or their inverses and these morphisms are natural by the very definition of the arrows in the categories $\tilde{R}$ and $\tilde{S}$. 

Let us now define a natural transformation $\eta:1_{\tilde{S}} \to \tilde{Z}\circ \tilde{W}$, which will serve as the unit of our adjunction. 

First, recall that if $X$ is an object of $\tilde{S}$ then the arrow 
\[
\eta'_{X}:X=(\pi_{S}^{{\cal H}}(X), \pi_{S}^{{\cal K}}(X)) \to \tilde{Z}(\tilde{W}(X))=(\pi_{S}^{{\cal H}}(\tilde{Z}(\tilde{W}(X))), \pi_{h}^{{\cal K}}(\tilde{Z}(\tilde{W}(X))))
\]
has as second component $\pi_{S}^{{\cal K}}(X)\to \pi_{S}^{{\cal K}}(\tilde{Z}(\tilde{W}(X)))$ the inverse of the isomorphism $\epsilon'_{Z(X)}:\pi_{R}^{{\cal K}}(\tilde{Z}(\tilde{W}(X))) \to \pi_{S}^{{\cal K}}(X)$. 

For any object $({\cal C}, {\cal D})$ in $\tilde{S}$, we define $\eta_{({\cal C}, {\cal D})}$ to be the triple
\[
\eta_{({\cal C}, {\cal D})}:=({\epsilon'_{W({\cal C}, {\cal D})}}^{-1}, \eta'_{({\cal C}, {\cal D})}, \eta'_{Z(W({\cal C}, {\cal D}))}).
\]
Let us check that this triple actually defines an arrow 
\[
({\cal C}, {\cal D}) \to \tilde{Z}\circ \tilde{W}({\cal C}, {\cal D})=(\pi_{S}^{\cal H}(Z(W({\cal C}, {\cal D}))), \pi^{{\cal K}}_{R}(W({\cal C}, {\cal D}))) 
\]
in the category $\tilde{S}$. This amounts to verifying that the both the squares in the following diagram are commutative. Note that, since $Z$ preserves the second component, 

\[
\begin{array}{lll}
W(\pi_{S}^{\cal H}(Z(W({\cal C}, {\cal D}))), \pi^{{\cal K}}_{R}(W({\cal C}, {\cal D}))) & = & W(\pi_{S}^{\cal H}(Z(W({\cal C}, {\cal D})), \pi_{S}^{\cal K}(Z(W({\cal C}, {\cal D}))) \\
& = & W(Z)({\cal C}, {\cal D}) 
\end{array}
\]
and hence in the diagram below we have 
\[
\pi^{{\cal K}}_{R}(W(\pi_{S}^{\cal H}(Z(W({\cal C}, {\cal D}))), \pi^{{\cal K}}_{R}(W({\cal C}, {\cal D}))))=\pi^{{\cal K}}_{R}(W(Z(W({\cal C}, {\cal D})))).
\]

\[  
\xymatrix {
J({\cal D}) \ar[d]^{\chi_{\cal D}} \ar[rrr]^{J(\eta'_{({\cal C}, {\cal D})})} & & & J(\pi^{{\cal K}}_{R}(W({\cal C}, {\cal D}))) \ar[d]^{\chi_{{\pi^{{\cal K}}_{R}(W({\cal C}, {\cal D}))}}}\\
I({\cal C}) \ar[d]^{\xi_{W({\cal C}, {\cal D})}} \ar[rrr]^{I({\epsilon'_{W({\cal C}, {\cal D})}}^{-1})} & & & I(\pi_{S}^{\cal H}(Z(W({\cal C}, {\cal D}))))  \ar[d]^{\xi_{\pi_{S}^{\cal H}(Z(W({\cal C}, {\cal D})))}}\\
J(\pi^{{\cal K}}_{R}(W(({\cal C}, {\cal D})))) \ar[rrr]^{J(\eta'_{Z(W({\cal C}, {\cal D}))})} & & & J(\pi^{{\cal K}}_{R}(W(Z(W({\cal C}, {\cal D})))))}.
\]
   
Now, the bottom square obviously commutes (cf. the diagram above made of two triangles, which we observed to commute), while the commutativity of the top square can be proved by observing that in the following diagram all the internal trapezoids are commutative whence the outer square is.

\[  
\xymatrix {
J({\cal D}) \ar[dr]^{\chi_{({\cal C}, {\cal D})}} \ar[ddd]^{\chi_{({\cal C}, {\cal D})}} \ar[rr]^{J(\eta'_{({\cal C}, {\cal D})})} & & J(\pi^{{\cal K}}_{R}(W({\cal C}, {\cal D}))) \ar[ddd]^{\chi_{{\pi^{{\cal K}}_{R}(W({\cal C}, {\cal D}))}}}\\
& I({\cal C}) \ar[d]^{\xi_{W({\cal C}, {\cal D})}} \ar[ur]^{\xi_{W({\cal C}, {\cal D})}} & \\
& J(\pi_{R}^{{\cal K}}(W({\cal C}, {\cal D}))) \ar[dr]^{\chi_{\pi_{R}^{{\cal K}}(W({\cal C}, {\cal D}))}} & \\
I(\pi^{{\cal K}}_{R}(W(({\cal C}, {\cal D})))) \ar[ur]^{\xi_{W({\cal C}, {\cal D})}} 
\ar[rr]_{I({\epsilon'_{W({\cal C}, {\cal D})}}^{-1})} & & I(\pi^{{\cal K}}_{R}(W(Z(W({\cal C}, {\cal D})))))}.
\]      

The naturality of $\eta$ is immediate from the fact that each of its components is a composition of (morphisms naturally isomorphic to) the morphisms $\xi$ and $\chi$ or their inverses, and these morphisms are natural by the very definition of the arrows in the categories $\tilde{R}$ and $\tilde{S}$.

So far, we have shown that $\eta$ and $\epsilon$ are natural transformations respectively $1_{\tilde{S}}\to \tilde{Z}\circ \tilde{W}$ and $\tilde{W}\circ \tilde{Z}\to 1_{\tilde{R}}$. To conclude that they are the unit and counit of an adjunction between the functors $\tilde{Z}$ and $\tilde{W}$, it remains to show that they satisfy the triangular identities. 

First, let us prove that for every $({\cal C}, {\cal D})\in \tilde{R}$ $\tilde{Z}(\epsilon_{({\cal C}, {\cal D})}) \circ \eta_{\tilde{Z}({\cal C}, {\cal D})}=1_{Z(({\cal C}, {\cal D}))}$. By definition of the functor $\tilde{Z}$, $\tilde{Z}(\epsilon_{({\cal C}, {\cal D})})=(\epsilon'_{W(Z(({\cal C}, {\cal D})))}, {\eta'_{Z({\cal C}, {\cal D})}}^{-1}, w)$, where $w$ is the unique (iso)morphism making the following diagram commute;

\[  
\xymatrix {
\pi^{\cal K}_{S}(Z(W(Z(({\cal C}, {\cal D}))))) \ar[d]^{\eta'_{Z(W(Z({\cal C}, {\cal D})))}} \ar[r]^{{\eta'_{Z({\cal C}, {\cal D})}}^{-1}} & {\cal D}' \ar[d]^{{\eta'_{Z({\cal C}, {\cal D})}}}\\
\pi^{{\cal K}}_{R}(W(Z(W(Z(({\cal C}, {\cal D})))))) \ar[r]^{w}   & \pi^{\cal K}_{S}(Z(W(Z(({\cal C}, {\cal D})))))   }
\]   
That is, $w={\eta'_{Z(W(Z({\cal C}, {\cal D})))}}^{-1}$.

Summarizing, $\tilde{Z}(\epsilon_{({\cal C}, {\cal D})})=(\epsilon'_{W(Z(({\cal C}, {\cal D})))}, {\eta'_{Z({\cal C}, {\cal D})}}^{-1}, {\eta'_{Z(W(Z({\cal C}, {\cal D})))}}^{-1})$.

On the other hand, we have 
\[
\eta_{Z(({\cal C}, {\cal D}))}:=({\epsilon'_{W(Z(({\cal C}, {\cal D})))}}^{-1}, \eta'_{Z(({\cal C}, {\cal D}))}, \eta'_{Z(W(Z(({\cal C}, {\cal D}))))}).
\]

From these two expressions we thus conclude that for any $({\cal C}, {\cal D})\in \tilde{R}$, $\tilde{Z}(\epsilon_{({\cal C}, {\cal D})}) \circ \eta_{\tilde{Z}({\cal C}, {\cal D})}=1_{Z(({\cal C}, {\cal D}))}$, as required. 

It remains to prove that the other triangular identity holds, i.e. that for any object $({\cal C}, {\cal D})$ of $\tilde{S}$, $\epsilon_{\tilde{W}({\cal C}, {\cal D})} \circ \tilde{W}(\eta_{({\cal C}, {\cal D})})=1_{W({\cal C}, {\cal D})}$. Now, we have that 

\[
\epsilon_{\tilde{W}({\cal C}, {\cal D})}=({\epsilon'_{W(({\cal C}, {\cal D}))}}, {\eta'_{Z(W({\cal C}, {\cal D}))}}^{-1}, {\epsilon'_{W(Z(W({\cal C}, {\cal D})))}})
\]

On the other hand, 

\[
\eta_{({\cal C}, {\cal D})}=({\epsilon'_{W({\cal C}, {\cal D})}}^{-1}, \eta'_{({\cal C}, {\cal D})}, \eta'_{Z(W({\cal C}, {\cal D}))}),
\]
whence $\tilde{W}(\eta_{({\cal C}, {\cal D})})=({\epsilon'_{W({\cal C}, {\cal D})}}^{-1}, \eta'_{Z(W({\cal C}, {\cal D}))})$.

From this it is immediate to conclude that $\epsilon_{\tilde{W}({\cal C}, {\cal D})} \circ \tilde{W}(\eta_{({\cal C}, {\cal D})})=1_{W({\cal C}, {\cal D})}$, as required. 

To finish the proof of the theorem, it remains to observe that the adjunction just generated is a reflection because $\eta_{({\cal C}, {\cal D})}$ is an isomorphism in $\tilde{R}$ for every object $({\cal C}, {\cal D})$ of $\tilde{R}$, since both of its components are isomorphisms.  

\end{proofs}

One might naturally wonder whether our method for building reflections (or coreflections) can be generalized to a method for building arbitrary adjunctions; the answer to this question is positive, but we shall not address this issue here since the details of such generalizations would bring us too far from the scope of this paper, and would be considerably more technical than the treatment just given for reflections, to the extent that the ratio between technical sophistication and conceptual understanding could be judged rather unbalanced by most readers. On the other hand, considering that every adjunction can be naturally obtained as a composite of a reflection with a coreflection (through the comma category construction), the choice of concentrating the analysis on reflections does not appear as a particularly restrictive one.

\subsection{Duality principles}\label{duality}

We note that, given a set of data 
\[
({\cal H}, {\cal K}, {\cal U}, I, J, R, S, Z, W)
\]
satisfying the conditions of our method, the set of data 
\[
({\cal K}, {\cal H}, {\cal U}, J, I, R^{\textrm{op}}, Z^{\textrm{op}}, S^{\textrm{op}}, W^{\textrm{op}}, Z^{\textrm{op}})
\]
also satisfies them.

There is another duality principle implicit in the context of our main result, which can be fruitfully exploited as a way of building reflections as instances of the theorem of the last section. We can illustrate this as follows.

Suppose having two categories $\cal H$ and $\cal K$, a category $\cal U$, two functors $I:{\cal H}\to {\cal U}$, $J:{\cal K}\to {\cal U}$ and two binary relations $R$ and $S$ on $Ob({\cal H})\times Ob({\cal K})$. In addition to this, suppose having, for every $({\cal C}, {\cal D})\in R$, an arrow $\alpha_{({\cal C}, {\cal D})}:  J({\cal C}, {\cal D})\to I({\cal C}, {\cal D})$ in $\cal U$, for every $({\cal C}, {\cal D})\in S$ an arrow $\beta_{({\cal C}, {\cal D})}: I({\cal C}, {\cal D}) \to J({\cal C}, {\cal D})$ in $\cal U$, and having furthermore two functions $Z:R\to S$ and $W:S\to R$ such that $Z$ keeps the second component fixed and $W$ keeps the first component fixed, which satisfy the following conditions:

\begin{enumerate}
\item For any $({\cal C}, {\cal D})\in R$, $Z(({\cal C}, {\cal D}))=(\pi^{{\cal H}}_{S}(Z(({\cal C}, {\cal D}))), {\cal D})$; we require the composite

\[
\alpha_{({\cal C}, {\cal D})}\circ \beta_{W_{R}({\cal C}, {\cal D})}: J(\pi^{{\cal H}}_{S}(W_{R}(({\cal C}, {\cal D})))) \to I({\cal C})
\]
to be induced by a (canonically chosen) morphism in $\cal H$ ${\lambda'_{({\cal C}, {\cal D})}}^{-1}: {\cal C} \to \pi^{{\cal H}}_{S}(Z(({\cal C}, {\cal D}))))$. 

\item For any $({\cal C}, {\cal D})\in S$, $W(({\cal C}, {\cal D}))=({\cal C}, \pi^{{\cal K}}_{R}(W(({\cal C}, {\cal D}))))$; we require the composite
\[
\beta_{({\cal C}, {\cal D})} \circ \alpha_{W({\cal C}, {\cal D})}: I(\pi^{{\cal K}}_{R}(W_{S}(({\cal C}, {\cal D})))) \to J({\cal D})
\]
to be induced by a (canonically chosen) morphism in $\cal K$ ${\mu'_{({\cal C}, {\cal D})}}: {\cal D}\to \pi^{{\cal K}}_{R}(W_{S}(({\cal C}, {\cal D})))$;

\item for any $({\cal C}, {\cal D})\in R_{R}$, the arrow
\[
\alpha_{({\cal C}, {\cal D})}: J({\cal D}) \to I({\cal C})
\]
is an isomorphism. [Note that, since $\lambda'_{({\cal C}, {\cal D})}$ is an isomorphism, 
\[
\beta_{Z({\cal C}, {\cal D})}: I(\pi^{{\cal H}}_{S}(Z(({\cal C}, {\cal D})))) \to J({\cal D})
\]
is an isomorphism as well.]

\item For any $({\cal C}, {\cal D})\in R$, $\mu'_{Z(({\cal C}, {\cal D}))}$ is an isomorphism.

\end{enumerate} 

Then the set of data 
\[
({\cal H}^{\textrm{op}}, {\cal K}^{\textrm{op}}, {\cal U}^{\textrm{op}}, I^{\textrm{op}}, J^{\textrm{op}}, R, S, Z, W)
\]
satisfies the hypotheses of our main theorem. 

The theorem thus yields a coreflection between the category $\tilde{R}$ and the category $\tilde{S}$, given by functors $\tilde{Z}:\tilde{R}\to \tilde{S}$ and $\tilde{W}:\tilde{S}\to \tilde{R}$, equivalently a reflection between the category $\tilde{R}^{\textrm{op}}$ and the category $\tilde{S}^{\textrm{op}}$, given by functors $\tilde{Z}^{\textrm{op}}:\tilde{R}^{\textrm{op}}\to \tilde{S}^{\textrm{op}}$ and $\tilde{W}^{\textrm{op}}:\tilde{S}^{\textrm{op}}\to \tilde{R}^{\textrm{op}}$.

The categories and functors involved in this adjunction can be described explicitly as follows.

The objects of $\tilde{R}^{\textrm{op}}$ are the elements of $R$ while the arrows $({\cal C}, {\cal D})\to ({\cal C}', {\cal D}')$ are the pairs $(u,v)$, where 
\[
u:{\cal C}\to {\cal C}'
\]
and
\[
v:{\cal D}\to {\cal D}'
\]
are arrows respectively in the categories $\cal H$ and $\cal K$ such that the following diagram commutes:

\[  
\xymatrix {
I({\cal C}') \ar[r]^{I(u)} & I({\cal C}) \\
J({\cal D}') \ar[u]^{\alpha_{({\cal C}', {\cal D}')}} \ar[r]^{J(v)} & J({\cal D}) \ar[u]^{\alpha_{({\cal C}, {\cal D})}}}
\]   

We shall occasionally write $(u,v,z)$ for $(u,v)$, where $z$ is the arrow 
\[
\pi^{{\cal H}}_{S}(Z(({\cal C}, {\cal D}))) \to \pi^{{\cal H}}_{S}(Z(({\cal C}', {\cal D}')))
\]
in $\cal H$ given by the factorization of $u$ across the isomorphisms $\mu'_{({\cal C}, {\cal D})}$ and $\mu'_{({\cal C}', {\cal D}')}$. 

The composition of arrows in $\tilde{R}^{\textrm{op}}$ is defined as the composition of the functors forming the various components.

Similarly, we define the category $\tilde{S}^{\textrm{op}}$. The objects of $\tilde{S}^{\textrm{op}}$ are the elements of $S$ while the arrows $({\cal C}, {\cal D})\to ({\cal C}' \to {\cal D}')$ are the triples $(z,v,w)$, where 
\[
v:{\cal D}\to {\cal D}',
\]
\[
z:{\cal C} \to {\cal C}'
\]
and 
\[
w:\pi^{{\cal K}}_{R}(W(({\cal C}, {\cal D}))) \to \pi^{{\cal K}}_{R}(W(({\cal C}', {\cal D}')))
\]
are morphisms respectively in the categories $\cal K$, $\cal H$ and $\cal K$ such that the two squares in the following diagram commute:

\[  
\xymatrix {
J({\cal D}')  \ar[r]^{J(v)} & J({\cal D})  \\
I({\cal C}') \ar[u]^{\beta_{({\cal C}', {\cal D}')}}  \ar[r]^{I(z)} & I({\cal C}) \ar[u]^{\beta_{({\cal C}, {\cal D})}} \\
J(\pi^{{\cal K}}_{R}(W(({\cal C}', {\cal D}')))) \ar[r]^{J(w)} \ar[u]^{\alpha_{({\cal C}', {\cal D}')}} & J(\pi^{{\cal K}}_{R}(W(({\cal C}, {\cal D})))) \ar[u]^{\alpha_{({\cal C}, {\cal D})}}}.
\]   

The composition of arrows in $\tilde{S}^{\textrm{op}}$ is defined as the composition of the functors forming the various components.

The functors $\tilde{Z}^{\textrm{op}}:\tilde{R}^{\textrm{op}}\to \tilde{S}^{\textrm{op}}$ and $\tilde{W}^{\textrm{op}}:\tilde{S}^{\textrm{op}}\to \tilde{R}^{\textrm{op}}$ can be described as follows.

For any object $({\cal C}, {\cal D})$ of $\tilde{R}$, $\tilde{Z}^{\textrm{op}}(({\cal C}, {\cal D}))=Z(({\cal C}, {\cal D}))$ and for any arrow $(u,v,z):({\cal C}, {\cal D}) \to ({\cal C}', {\cal D}')$ in $\tilde{R}^{\textrm{op}}$ we set $\tilde{Z}^{\textrm{op}}((u,v,z))$ equal to the triple $(z, v, w)$, where $w:\pi^{{\cal K}}_{R}(W(Z({\cal C}, {\cal D})) \to \pi^{{\cal K}}_{R}(W(Z({\cal C}', {\cal D})')$ is the only arrow in $\cal K$ making the following diagram commute. 

\[  
\xymatrix {
\pi^{{\cal K}}_{R}(W(Z({\cal C}, {\cal D}))) \ar[r]^{w}  & \pi^{{\cal K}}_{R}(W(Z({\cal C}', {\cal D})')) \\
{\cal D} \ar[u]^{\lambda'_{Z({\cal C}, {\cal D})}} \ar[r]^{v} & {\cal D}' \ar[u]^{\lambda'_{Z({\cal C}', {\cal D}')}}  }
\]   

For any $({\cal C}, {\cal D})$ in $\tilde{S}$, $\tilde{W}^{\textrm{op}}(({\cal C}, {\cal D}))=W(({\cal C}, {\cal D}))$ and for any arrow $(z,v,w):({\cal C}, {\cal D})\to ({\cal C}', {\cal D}')$ in $\tilde{S}$, $\tilde{W}((z,v,w))=(z,w)$.

\subsection{From the relational to the functional context}\label{functional} 

In this section we specialize our general method for generating reflections to the context of functional relations $R$ and $S$. In fact, the main reason behind our choice of a relational context for formulating our theorem in section \ref{general} is the fact that this more general context provides us with duality principles (cf. section \ref{duality} above) which do not hold in the restricted functional context.

Suppose having two categories $\cal H$ and $\cal K$, a category $\cal U$, two functors $I:{\cal H}\to {\cal U}$, $J:{\cal K}\to {\cal U}$, two functions $f:Ob({\cal H})\to Ob({\cal K})$, $g:Ob({\cal K})\to Ob({\cal H})$, for every object ${\cal C}\in {\cal H}$ an arrow 
\[
\xi_{{\cal C}}:I({\cal C})\to J(f({\cal C}))
\]
in $\cal U$, and for every element ${\cal D}\in {\cal K}$ an arrow 
\[
\chi_{{\cal D}}:J({\cal D})\to I(g({\cal D}))
\]
in $\cal U$.

Assume that the composite $\chi_{f({\cal C})}\circ \xi_{{\cal C}}$ is of the form $I({\eta_{\cal C}})$ for some morphism $\eta_{\cal C}:{\cal C} \to g(f({\cal C}))$, while the composite $\xi_{g({\cal D})}\circ \chi_{{\cal D}}:J({\cal D})\to I(f(g({\cal D})))$ is of the form $J({\epsilon_{{\cal D}}}^{-1})$ for an isomorphism $\epsilon_{{\cal D}}:f(g({\cal D})) \to {\cal D}$, as in the following diagrams:
\[  
\xymatrix {
I({\cal C}) \ar[r]^{\xi_{{\cal C}}} \ar[dr]^{I({\eta'_{\cal C}})} & J(f({\cal C})) \ar[d]^{\chi_{f({\cal C})}}\\
& I(g(f({\cal C})))}
\]

\[  
\xymatrix {
J({\cal D}) \ar[r]^{\chi_{{\cal D}}} \ar[dr]^{J({\epsilon'_{\cal D}}^{-1})} & I(g({\cal D})) \ar[d]^{\xi_{g({\cal D})}}\\
& \Sh(f(g({\cal D})))}
\]

Suppose moreover that $\xi_{{\cal C}}$ is an isomorphism for every $\cal C$; note that this implies, by definition of $\eta_{\cal C}$, that $\chi_{f({\cal C})}$ is an isomorphism for each ${\cal C}$ in $\cal H$. 

Out of these data, we can construct two sets $R_{R}$ and $R_{S}$ and two functions $Z_{R}:R_{R}\to R_{S}$ and $Z_{S}:R_{S}\to R_{R}$ satisfying the hypotheses of our main theorem, as follows.

\begin{enumerate}

\item We define the relation $R$ as the graph $R_{f}$ of a function $f:Ob({\cal H})\to Ob({\cal K})$, i.e. the set of pairs of the form $({\cal C}, f({\cal C}))$ for ${\cal C}\in {\cal H}$, and $S$ as the inverse of the graph $R_{g}$ of a function $g:Ob({\cal K})\to Ob({\cal H})$, i.e. the set of pairs of the form $(g({\cal D}), {\cal D})$ for ${\cal D}\in {\cal K}$. 

\item We define $Z:R_{f}\to R_{g}$ as the function sending a pair $({\cal C}, {\cal D})$ in $R_{f}$ to the pair $(g({\cal D}), {\cal D})$.

\item We define $W:R_{g}\to R_{f}$ as the function sending a pair $({\cal C}, {\cal D})$ in $R_{g}$ to the pair $({\cal C}, f({\cal C}))$.
\end{enumerate}

These data satisfy the hypotheses of the theorem with respect to the morphisms $\xi_{({\cal C}, {\cal D})}$ and $\chi_{({\cal C}, {\cal D})}$ defined by:

\begin{enumerate}
\item For any $({\cal C}, {\cal D})\in R_{f}$, $\xi_{({\cal C}, {\cal D})}=\xi_{{\cal C}}$;

\item For any $({\cal C}, {\cal D})\in R_{g}$, $\chi_{({\cal C}, {\cal D})}=\chi_{{\cal D}}$.
\end{enumerate}

Our general theorem of section \ref{general} thus yields a coreflection between the category $\tilde{R_{f}}$ and the category $\tilde{R_{g}}$. Under the following additional assumptions, the description of the categories $\tilde{R_{f}}$ and $\tilde{R_{g}}$ and of the functors $\tilde{Z_{f}}:\tilde{R_{f}}\to \tilde{R_{g}}$ radically simplifies. We assume the following conditions:

\begin{enumerate}

\item For any morphism $x:{\cal C} \to {\cal C}'$ in $\cal H$ there is at most one morphism $y:f({\cal C}) \to f({\cal C}')$ in $\cal K$ such that the diagram 

\[  
\xymatrix {
I({\cal C}) \ar[d]^{\xi_{{\cal C}}} \ar[r]^{I(x)} & I({\cal C}')  \ar[d]^{\xi_{{\cal C}'}}\\
J(f({\cal C})) \ar[r]_{J(y)} & J(f({\cal C}))}
\]  
commutes;

\item For any morphism $y:{\cal D} \to {\cal D}'$ in $\cal K$ there is at most one morphism $x:g({\cal D}) \to g({\cal D}')$ in $\cal H$ such that the diagram 

\[  
\xymatrix {
J({\cal D}) \ar[d]^{\chi_{{\cal D}}} \ar[r]^{J(y)} & J({\cal D}')  \ar[d]^{\chi_{{\cal D}'}}\\
I(g({\cal D})) \ar[r]_{I(x)} & I(g({\cal D}'))}
\]  
commutes. 
\end{enumerate} 

Under these hypotheses, the category $\tilde{R}$ is equivalent to the category $\tilde{\cal H}$ whose objects are the objects of $\cal H$ and whose arrows ${\cal C}\to {\cal C}'$ are the arrows $s:{\cal C}\to {\cal C}'$ in $\cal H$ such that there is a unique arrow $t_{s}:f({\cal C}) \to f({\cal C}')$ in $\cal K$ such that the diagram

\[  
\xymatrix {
I({\cal C}) \ar[d]^{\xi_{({\cal C}, {\cal D})}} \ar[r]^{I(s)} & I({\cal C}')  \ar[d]^{\xi_{({\cal C}', {\cal D}')}}\\
J(f({\cal C})) \ar[r]_{J(t_{s})} & J(f({\cal C}'))}
\]  

commutes.

Under the same hypotheses, the category $\tilde{S}$ can be identified with the category $\tilde{\cal K}$ having as objects the objects of $\cal K$ and as arrows ${\cal D}\to {\cal D}'$ the arrows $t:{\cal D} \to {\cal D}'$ in ${\cal K}$ such that there is a (unique) arrow $r_{t}:g({\cal D}) \to g({\cal D}')$ making the diagram

\[  
\xymatrix {
J({\cal D}) \ar[d]^{\chi_{({\cal C}, {\cal D})}} \ar[r]^{J(t)} & J({\cal D}')  \ar[d]^{\chi_{({\cal C}', {\cal D}')}}\\
I(g({\cal D})) \ar[r]_{I(r_{t})} & I(g({\cal D}'))}
\]  

commute, and the morphism $r_{t}$ defines an arrow $g({\cal D}) \to g({\cal D}')$ in $\tilde{{\cal H}}$, equivalently there exists a unique arrow 
\[
y:f(g({\cal D})) \to f(g({\cal D}'))
\]
such that the diagram

\[  
\xymatrix {
I(g({\cal D})) \ar[d]^{\xi_{g({\cal D})}} \ar[r]^{I(r_{t})} & I(g({\cal D}'))  \ar[d]^{\xi_{g({\cal D}')}}\\
J(f(g({\cal D}))) \ar[r]_{J(y)} & J(f(g({\cal D}')))}
\]  
commutes. 

The functor $\tilde{Z_{f}}:{\cal H}\to {\cal K}$ can be characterized as the functor sending any arrow $s$ in $\cal H$ to the arrow $t_{s}$ defined above, while the functor $\tilde{Z_{g}}:{\cal K}\to {\cal H}$ sends any arrow $r$ in $\cal K$ to the arrow $r_{t}$ defined above. 

Notice that, by the remarks of section \ref{duality}, a dual version of this corollary yielding reflections (in the functional context) also holds.

\subsection{Completeness}\label{completeness}

In this section we show that every reflection between categories can be obtained as an application of our general method. 

Let $\cal A$ and $\cal B$ categories and $R:{\cal A}\to {\cal B}$, $L:{\cal B}\to {\cal A}$ be two adjoint functors ($L$ left adjoint and $R$ right adjoint) with unit $\eta:1_{{\cal B}}\to R\circ L$ and counit $\epsilon:L\circ R\to 1_{{\cal A}}$. If this adjunction is a coreflection then we can obtain it as the result of applying our main theorem to the following set of data.

We define ${\cal H}={\cal A}$, ${\cal K}={\cal B}$, ${\cal U}={\cal B}$, $I$ as the functor $R:{\cal H}={\cal A}\to {\cal B}={\cal U}$, $J$ as the functor $1_{{\cal B}}:{\cal K}={\cal B}\to {\cal B}={\cal U}$, $f$ as the underlying function on objects of the functor $R$, $g$ as the underlying function on objects of the functor $L$, $\xi_{a}$ (for $a\in {\cal A}$) as the identity arrow $1_{R(a)}:I(a)=R(a)\to R(a)=J(f(a))$, $\chi_{b}$ (for $b\in {\cal B}$) as the arrow $\eta_{b}:J(b)=b\to R(L(b))=I(g(b))$. 

It is immediate to verify that this set of data satisfies the hypotheses of our method (rewritten in the functional context as in section \ref{functional}). The resulting category $\tilde{R_{f}}$ is clearly isomorphic to $\cal A$ (since the $\xi$ are all identities), while the category $\tilde{R_{g}}$ is isomorphic to $\cal B$, since for any arrow $u:b\to b'$ in $\cal B$ there exists exactly one arrow $v:g(b)\to g(b')$ such that the following diagram commutes, namely $L(u)$.
    
\[  
\xymatrix {
b \ar[d]^{\chi_{b}} \ar[r]^{J(u)} & b' \ar[d]^{\chi_{b'}}\\
I(g(b)) \ar[r]^{I(v)} & I(g(b')) }
\]     

Indeed, this diagram is precisely the naturality square for the unit $\eta$ with respect to the arrow $u$, if $v=L(u)$.

Another way of recovering this adjunction as an application of our method consists in selecting a different set of data leading to the same reflection: one can alternatively choose ${\cal H}={\cal A}$, ${\cal K}={\cal B}$, ${\cal U}={\cal A}$, $I$ as the functor $1_{\cal A}:{\cal H}={\cal A}\to {\cal A}={\cal U}$, $J$ as the functor $L:{\cal K}={\cal B}\to {\cal A}={\cal U}$, $f$ as the underlying function on objects of the functor $R$, $g$ as the underlying function on objects of the functor $L$, $\xi_{a}$ (for $a\in {\cal A}$) as the arrow ${\epsilon_{a}}^{-1}:I(a)=a\to L(R(a))=J(f(a))$, $\chi_{b}$ (for $b\in {\cal B}$) as the identity arrow arrow $1_{L(b)}:J(b)=L(b)\to L(b)=I(g(b))$.

\section{Reflections from geometric morphisms}\label{geomrefl}

In this section we apply our general method for generating reflections to a specific context, namely the topos-theoretic interpretation of Stone-type dualities established in \cite{OC11}. We shall be able to naturally recover as special applications of our method all the Stone-type reflections or coreflections that we discussed in \cite{OC11}, providing at the same time a uniform way for building such adjunctions.

\begin{enumerate}

\item Suppose that $\cal H$ is a category of poset structures ${\cal C}$, whose arrows are precisely the monotone maps ${\cal C}\to {\cal C}'$, and that $\cal K$ is a category of poset structures ${\cal D}$, each of which equipped with a subcanonical Grothendieck topology $K_{{\cal D}}$, whose morphisms ${\cal D}\to {\cal D}'$ are precisely the morphisms of sites $({\cal D}, K_{\cal D})\to ({\cal D}', K_{{\cal D}'})$. If we take $\cal U$ to be the (skeleton of the) category of Grothendieck toposes then we have two functors 
\[
I:{\cal H}\to {\cal U}
\]
and
\[
J:{\cal K}^{\textrm{op}}\to {\cal U},
\]
defined as follows.

\begin{enumerate}
\item The functor $I:{\cal K} \to {\cal U}$ sends a poset $\cal C$ in $\cal K$ to the topos $[{\cal C}, \Set]$ and an arrow $v:{\cal C}\to {\cal C}'$ in $\cal K$ to the geometric morphism 
\[
E(v):[{\cal C}, \Set]\to [{\cal C}', \Set]
\]
canonically induced by $v$.

\item The functor $J:{\cal K}^{\textrm{op}}\to {\cal U}$ sends a category ${\cal D}$ in $\cal K$ to the topos $\Sh({\cal D}, K_{\cal D})$ and a morphisms of sites $u:({\cal D}, K_{\cal D})\to ({\cal D}', K_{{\cal D}'})$ to the induced geometric morphism 
\[
\Sh(u):\Sh({\cal D}', K_{{\cal D}'}) \to \Sh({\cal D}, K_{\cal D}).
\] 
\end{enumerate}
If one has two functions $f:Ob({\cal H})\to Ob({\cal K})$ and $g:Ob({\cal K})\to Ob({\cal H})$ and geometric morphisms 
\[
\alpha_{\cal C}:\Sh(f({\cal C}), J_{f({\cal C})}) \to [{\cal C}, \Set]
\]
(for $\cal C$ in $\cal H$) and 
\[
\beta_{\cal D}:[g({\cal D}), \Set] \to \Sh({\cal D}, J_{{\cal D}})
\]
(for $\cal D$ in $\cal K$) satisfying the hypotheses of theorem of section \ref{functional} then there is a coreflection between ${\cal H}$ and ${\cal K}^{\textrm{op}}$ given by functors $\tilde{Z}:\tilde{R_{f}}\to \tilde{R_{g}}$ and $\tilde{W}:\tilde{R_{g}}\to \tilde{R_{f}}$. As shown in the previous section, under the following hypotheses the categories $\tilde{R_{f}}$ and $\tilde{R_{g}}$ admit simpler descriptions: 

\begin{enumerate}

\item For any monotone map $y:{\cal C} \to {\cal C}'$ there is at most one morphism of sites $x:(f({\cal C})', K_{f({\cal C}')}) \to (f({\cal C}, K_{f({\cal C})}))$ such that the diagram 

\[  
\xymatrix {
\Sh(f({\cal C}), K_{f({\cal C})}) \ar[d]^{\alpha_{{\cal C}}}  \ar[r]_{\Sh(x)} & \Sh(f({\cal C}'), K_{f({\cal C}')}) \ar[d]^{\alpha_{{\cal C}'}}\\
[{\cal C}, \Set] \ar[r]^{E(y)} & [{\cal C}', \Set] }
\]  
commutes; 

\item For any morphism of sites $x:({\cal D}, K_{\cal D}) \to ({\cal D}', K_{{\cal D}'})$ there is at most one monotone map $y:g({\cal D})' \to g({\cal D})$ such that the diagram 

\[  
\xymatrix {
[g({\cal D}'), \Set] \ar[d]^{\beta_{{\cal D}'}} \ar[r]_{E(y)} & \ar[d]^{\beta_{{\cal D}}} [g({\cal D}), \Set]\\
\Sh({\cal D}', K_{{\cal D}'})  \ar[r]^{\Sh(x)} & \Sh({\cal D}, K_{{\cal D}})}
\]  
commutes;
\end{enumerate} 

We note that condition $(a)$ is always satisfied, the $\alpha$ are being equivalences. Indeed, if $\cal A$ and $\cal B$ are Cauchy-complete categories (e.g. posets), a functor $l:{\cal A}\to {\cal B}$ can be recovered up to isomorphism from the associated geometric morphism $E(l):[{\cal A}, \Set]\to [{\cal B}\to \Set]$ as the restriction to the full subcategories spanned by the representables of the left adjoint to the inverse image functor of $E(l)$. 

We can identify a large class of naturally arising contexts in which condition $(b)$ is satisfied. Suppose that for every $\cal D$ in $\cal K$, $g({\cal D})$ is a subset of ${\cal D}$ such that the map ${\cal D}\to Id(g({\cal D}))$ sending any element $d\in {\cal D}$ to the ideal $\{d'\in g({\cal D}) \textrm{ | } d'\leq d\}$ is a flat $J_{g({\cal D})}$-continuous flat functor $F_{\cal D}$ corresponding via Diaconescu's equivalence to the geometric morphism $\beta_{\cal D}:[{\cal D}, \Set]\to \Sh(g({\cal D}), J_{g({\cal D})})$; then condition $(b)$ is satisfied. Indeed the commutativity of the square in condition $(b)$ is equivalent to the commutativity of the following diagram
\[  
\xymatrix { 
g({\cal D}') \ar[r]^{y} \ar[d]^{F_{{\cal D}'}} & g({\cal D}) \ar[d]^{F_{{\cal D}}}\\
Id({\cal D}') \ar[r]^{Id(x)} & Id({\cal D}) }
\]  
and $Id(x)$, being a frame homomorphism, is forced to be the function sending an ideal $I$ in ${\cal D}'$ to the union of the ideals in ${\cal D}$ of the form $F_{{\cal D}}(x(s))$ for $s$ such that $F_{{\cal D}'}(s)\subseteq I$. As we already remarked above, $x$ is uniquely determined by $Id(x)$.

Under these hypotheses, if the Grothendieck topologies $K_{\cal D}$ are $C$-induced for an invariant $C$ satisfying the conditions in the theorems of \cite{OC11} then $\tilde{R_{f}}$ can be described as the category whose objects are the posets in $\cal H$ and whose morphisms ${\cal C}\to {\cal C}'$ are the morphisms in $\cal H$ whose corresponding frame homomorphisms $Id({\cal C}')\to Id({\cal C})$ send $C$-compact elements to $C$-compact elements, while $\tilde{R_{g}}$ can be described as the category whose objects are the posets in $\cal K$ and whose arrows are the arrows ${\cal D}\to {\cal D}'$ in $\cal K$ such that the map $Id(g({\cal D}))\to Id(g({\cal C}'))$ sending an ideal $I$ on $g({\cal D})$ to the union of the ideals of the form $I_{s}:=\{x\in g({\cal D}') \textrm{ | } x\leq g(s)\}$ for $s\in I$ is a complete frame homomorphism.
  
The two adjunctions for atomic frames and locally connected frames obtained in \cite{OC11}, as well as the Lindenbaum-Tarski adjunction between sets and complete Boolean algebras, are particular instances of this kind of adjunctions.

\item Suppose that $\cal H$ is a category of poset structures ${\cal C}$, each of which equipped with a subcanonical Grothendieck topology $J_{{\cal C}}$, whose morphisms ${\cal C}\to {\cal C}'$ are precisely the morphisms of sites $({\cal C}, J_{\cal C})\to ({\cal C}', J_{{\cal C}'})$, and that $\cal K$ is a category of poset structures ${\cal D}$, each of which equipped with a subcanonical Grothendieck topology $K_{{\cal D}}$, whose morphisms ${\cal D}\to {\cal D}'$ are precisely the morphisms of sites $({\cal D}, K_{\cal D})\to ({\cal D}', K_{{\cal D}'})$. If we take $\cal U$ to be the (skeleton of the) category of Grothendieck toposes then we have two functors 
\[
I:{\cal H}^{\textrm{op}}\to {\cal U}
\]
and
\[
J:{\cal K}^{\textrm{op}}\to {\cal U},
\]
defined as follows.

\begin{enumerate}
\item The functor $I:{\cal H}^{\textrm{op}}\to {\cal U}$ sends a category ${\cal C}$ in $\cal H$ to the topos $\Sh({\cal C}, J_{\cal C})$ and a morphisms of sites $s:({\cal C}, J_{\cal C})\to ({\cal C}', J_{{\cal C}'})$ to the induced geometric morphism 
\[
\Sh(s):\Sh({\cal C}', J_{{\cal C}'}) \to \Sh({\cal C}, K_{\cal C}).
\]

\item The functor $J:{\cal K}^{\textrm{op}}\to {\cal U}$ sends a category ${\cal D}$ in $\cal K$ to the topos $\Sh({\cal D}, K_{\cal D})$ and a morphisms of sites $u:({\cal D}, K_{\cal D})\to ({\cal D}', K_{{\cal D}'})$ to the induced geometric morphism 
\[
\Sh(u):\Sh({\cal D}', K_{{\cal D}'}) \to \Sh({\cal D}, K_{\cal D}).
\] 
\end{enumerate}

Suppose that one has two functions $f:Ob({\cal H})\to Ob({\cal K})$ and $g:Ob({\cal K})\to Ob({\cal H})$ and geometric morphisms 
\[
\xi_{\cal C}:\Sh({\cal C}), J_{{\cal C}}) \to \Sh(f({\cal C}), K_{({\cal C})})
\]
(for $\cal C$ in $\cal H$) and 
\[
\chi_{\cal D}: \Sh({\cal D}, K_{{\cal D}}) \to \Sh(g({\cal D}), J_{g({\cal D})}) 
\]
(for $\cal D$ in $\cal K$) satisfying the conditions of our method (specialized to the functional context as in section \ref{functional}). Then there is a reflection between ${\cal H}^{\textrm{op}}$ and ${\cal K}^{\textrm{op}}$ (equivalently, a coreflection between the categories $\cal H$ and $\cal K$) given by functors $\tilde{Z}:\tilde{R_{f}}\to \tilde{R_{g}}$ and $\tilde{W}:\tilde{R_{g}}\to \tilde{R_{f}}$. As shown above, under the following hypotheses the categories $\tilde{R_{f}}$ and $\tilde{R_{g}}$ admit simpler descriptions: 

\begin{enumerate}

\item For any morphism of sites $y:({\cal C}, J_{\cal C}) \to ({\cal C}', J_{{\cal C}'})$ there is at most one morphism of sites $x:(f({\cal C}), K_{f({\cal C}')}) \to (f({\cal C}'), K_{f({\cal C})})$ such that the diagram 

\[  
\xymatrix {
\Sh({\cal C}', J_{{\cal C}'}) \ar[d]^{\xi_{{\cal C}'}} \ar[r]_{\Sh(y)} & \Sh({\cal C}, J_{\cal C}) ) \ar[d]^{\xi_{{\cal C}'}}\\
\Sh(f({\cal C}'), K_{f({\cal C}')}) \ar[r]^{\Sh(x)} & \Sh(f({\cal C}'), K_{f({\cal C}')} }
\]  
commutes; 

\item For any morphism of sites $x:({\cal D}, K_{\cal D}) \to ({\cal D}', K_{{\cal D}'})$ there is at most one morphism of sites $y:(g({\cal D}), J_{g({\cal D})})  \to (g({\cal D}', J_{g({\cal D}')})$ such that the diagram 

\[  
\xymatrix {
\Sh({\cal D}', K_{{\cal D}'}) \ar[d]^{\chi_{{\cal D}'}}  \ar[r]^{\Sh(x)} & \Sh({\cal D}, K_{{\cal D}}) \ar[d]^{\chi_{{\cal D}}}  \\
\Sh(g({\cal D}'), J_{g({\cal D}')})  \ar[r]_{\Sh(y)} & \Sh(g({\cal D}), J_{g({\cal D})})}
\]  
commutes.
\end{enumerate} 

We note that condition $(a)$ is always satisfied, the $\xi$ are being equivalences. Indeed, it is well-known that a morphism of subcanonical sites can be recovered up to isomorphism (in particular, uniquely, in case of morphisms between posets) from the corresponding geometric morphism as the restriction of its inverse image functor of the morphism to the full subcategories spanned by the representables. 

We can identify a large class of naturally arising contexts in which condition $(b)$ is satisfied. Suppose that for every $\cal D$ in $\cal K$, $g({\cal D})$ is a subset of $\cal D$ such that the inclusion $g({\cal D})\hookrightarrow {\cal D}$ yields a morphism of sites $(g({\cal D}), J_{g({\cal D})})\to ({\cal D}, K_{\cal D})$ (note that, for example, this condition always holds if $J_{g({\cal D})}$ is the Grothendieck topology induced by $K_{\cal D}$ on $g({\cal D})$) which induces the geometric morphism $\chi_{{\cal D}}$; then the commutativity of the diagram forces $y$ to be equal to the restriction of $x$ to $g({\cal D})$ and $g({\cal D}')$ and hence condition $(b)$ is satisfied.

Under these hypotheses, if the Grothendieck topologies $K_{\cal D}$ are $C$-induced for an invariant $C$ satisfying the conditions in the theorems of \cite{OC11} then $\tilde{R_{f}}$ can be described as the category whose objects are the posets in $\cal H$ and whose morphisms ${\cal C}\to {\cal C}'$ are the morphisms in $\cal H$ whose corresponding frame homomorphisms $Id_{J_{\cal C}}({\cal C})\to Id_{J_{{\cal C}'}}({\cal C}')$ send $C$-compact elements to $C$-compact elements, while $\tilde{R_{g}}$ can be described as the category whose objects are the posets in $\cal K$ and whose arrows ${\cal D}\to {\cal D}'$ in $\cal K$ are the morphism of sites $({\cal D}, K_{\cal D})\to ({\cal D}', K_{{\cal D}'})$ which restrict to a morphism of sites 
\[
(g({\cal D}), J_{g({\cal D})})\to (g({\cal D}'), J_{g({\cal D}')})
\]
(note that, in case for every $\cal D$ $J_{g({\cal D})}$ is the Grothendieck topology induced by $K_{\cal D}$ on $g({\cal D})$, for this condition to hold it suffices to require that the underlying function ${\cal D}\to {\cal D}'$ of the morphism restricts to a function $g({\cal D})\to g({\cal D}')$).
which is an arrow in $\cal H$.

Of course, another case in which condition $(b)$ is satisfied is when the $\xi_{\cal D}$ are equivalences (cf. the example below). 
 
As a useful illustration of this kind of adjunctions, we describe the following context. Let $\cal H$ be a full subcategory of the category of frames and $\cal K$ be a category of posets, each of which equipped with a subcanonical topology, whose arrows are the morphisms of the associated sites; we denote by $K_{\cal D}$ the Grothendieck topology associated to a poset $\cal D$ in $\cal K$. Given a frame ${\cal C}$ in $\cal H$, we denote by $J_{{\cal C}}$ the canonical topology on $\cal C$. Suppose that we have a function $f:Ob({\cal H})\to Ob({\cal K})$ with the property that for any $\cal C$ in $\cal H$, $f({\cal C})$ is a basis of $\cal C$ such that $K_{f({\cal C})}=J_{\cal C}|_{f({\cal C})}$, and denote by $g:Ob({\cal K})\to Ob({\cal H})$ the function sending a poset $\cal D$ in $\cal K$ to the frame $g({\cal D}):=Id_{K_{\cal C}}({\cal D})$. Suppose that for every $\cal D$ in $\cal K$ the canonical morphism ${\cal D} \to Id_{K_{\cal C}}({\cal D})$ factors through the inclusion $f(Id_{K_{\cal C}}({\cal D})) \hookrightarrow Id_{K_{\cal C}}({\cal D})$. Then the equivalences $\xi_{{\cal C}}:\Sh({\cal C}, J_{\cal C})\to \Sh(f({\cal C}), K_{f({\cal C})})$ and $\Sh({\cal D}, K_{\cal D})\simeq \Sh(Id_{K_{\cal C}}({\cal D}), J_{Id_{K_{\cal D}}({\cal D})})$ induced by the Comparison Lemma satisfy the hypotheses of our method. Since conditions $(a)$ and $(b)$ above are trivially satisfied, we obtain a reflection between the categories $\tilde{R_{f}}$ and $\tilde{R_{g}}$ given by the functors $\tilde{Z_{f}}$ and $\tilde{Z_{g}}$. The category $\tilde{R_{f}}$ has as objects the frames in $\cal H$ and as arrows the frame homomorphisms ${\cal C}\to {\cal C}'$ which factor through the inclusions $f({\cal C})\hookrightarrow {\cal C}$ and $f({\cal C}')\hookrightarrow {\cal C}'$, while the category $\tilde{R_{g}}$ has as objects the posets in $\cal K$ and as arrows ${\cal D}\to {\cal D}'$ the morphisms such that the corresponding morphism $Id_{K_{\cal D}}({\cal D})\to Id_{K_{\cal D}'}({\cal D}')$ factors through the inclusions $f(Id_{K_{\cal D}}({\cal D}))\hookrightarrow Id_{K_{\cal D}}({\cal D})$ and $f(Id_{K_{\cal D}'}({\cal D}'))\hookrightarrow Id_{K_{\cal D}'}({\cal D}')$. The functor $\tilde{Z_{f}}:\tilde{R_{f}}\to \tilde{R_{g}}$ sends a frame $\cal C$ to the poset $f({\cal C})$ and a frame homomorphism ${\cal C}\to {\cal C}'$ to its restriction $f({\cal C})\to f({\cal D})$, while the functor $\tilde{Z_{g}}:\tilde{R_{g}}\to \tilde{R_{f}}$ sends any arrow in $\cal K$ to the corresponding frame homomorphism. 

An example of a reflection of this kind is the reflection between the category of frames and the category of Boolean algebras providing a localic version of the Stone adjunction between the category of Boolean algebras and the opposite of the category of topological spaces (cf. \cite{OC11}); indeed, for any frame $F$, denoted by $F_{c}$ the Boolean algebra of its complemented elements, we have a geometric morphism $\Sh(F_{J_{f}})\to \Sh(F_{c}, J^{coh}_{F_{c}})$ induced by the inclusion $(F_{c}, J^{coh}_{F_{c}})\hookrightarrow (F, J_{F})$ (where $J_{F}$ is the canonical topology on $F$ and $J^{coh}_{F_{c}}$ is the coherent topology on $F_{c}$) and for any Boolean algebra $B$ we have an equivalence $\Sh(F_{B}, J_{F_{b}})\simeq \Sh(B, J^{coh}_{B})$.

\item

Our general method for building reflections can be profitably used also for establishing reflections with categories of topological spaces. For these purposes, it is often useful to select as category $\cal U$ the category of toposes paired with points defined in \cite{OC11}. Let us give a couple of examples of these kind of adjunctions.

\begin{enumerate}

\item The well-known \emph{Stone adjunction} between the category of Boolean algebras and the opposite of the category of topological spaces can be obtained as follows. Take $\cal H$ to be the opposite of the category of Boolean algebras and $\cal K$ to be the category of topological spaces. Define $\cal U$ to be the category of toposes paired with points (as defined in \cite{OC11}). We have two functors $I:{\cal H}\to {\cal U}$ and $J:{\cal K}\to {\cal U}$; the functor $I$ sends a Boolean algebra $B$ to the pair $(\Sh(B, J_{B}), p_{B})$ where $p_{B}$ is the set of all the points of the topos $\Sh(B, J_{B})$ and acts on the arrows accordingly, while the functor $J$ sends a topological space $X$ to the pair $(\Sh(X), P_{X})$ where $P_{X}$ is the set of points of $\Sh(X)$ indexed by the elements of $X$ (as in \cite{OC11}). The functions $f:Ob({\cal H})\to Ob({\cal K})$ and $g:Ob({\cal K})\to Ob({\cal H})$ can be defined as follows; $f$ sends a Boolean algebra $B$ to its Stone spectrum $X_{B}$, while $g$ sends a topological space to the Boolean algebra $X_{cl}$ of its clopen subsets. For each $B$ in $H$ we have an isomorphism $\xi_{B}:I(B)=(\Sh(B, J_{B}), p_{B}))\to (\Sh(X_{B}), P_{X_{B}})=J(f(B))$ in $\cal U$, while for any $X$ in $\cal K$ we have an arrow $J(X)=(\Sh(X), P_{X})\to (\Sh(X_{cl}, J_{X_{cl}}), p_{X_{cl}}) =I(g(X))$ in $\cal U$. It is immediate to see that this set of data satisfies the hypotheses of our method, from which we conclude that we have a reflection between ${\cal H}$ and $\cal K$, as required.     

\item The \emph{Alexandrov adjunction} between the category of preorders and the category of topological spaces can be obtained as follows. Take $\cal H$ to be the category of preorders, and $\cal K$ to be the category of topological spaces. Define $\cal U$ to be the category of toposes paired with points (as defined in \cite{OC11}). We have two functors $I:{\cal H}\to {\cal U}$ and $J:{\cal K}\to {\cal U}$ defined as follows. The functor $I$ sends a preorder $\cal P$ to the pair $([{\cal P}, \Set], E_{\cal P})$, where $E_{\cal P}$ is the indexing of points of $[{\cal P}, \Set]$ by elements of $\cal P$ (as in \cite{OC11}), acting on the arrows in the obvious way. The functor $J$ sends a topological space $X$ to the pair $(\Sh(X), P_{X})$ where $P_{X}$ is the set of points of $\Sh(X)$ indexed by the elements of $X$ (as in \cite{OC11}). The functions $f:Ob({\cal H})\to Ob({\cal K})$ and $g:Ob({\cal K})\to Ob({\cal H})$ are defined as follows; $f$ sends a preorder $\cal P$ to the Alexandrov space $A_{\cal P}$ based on $\cal P$ (i.e. the preorder $\cal P$ endowed with the Alexandrov topology), while the function $g$ sends a topological space $X$ to its specialization preorder $X_{\leq}$. For each $\cal P$ in $\cal H$ we have an isomorphism $\alpha_{\cal P}:(\Sh(A_{\cal P}), P_{A_{\cal P}})=J(f({\cal P}))\to I({\cal P})=[{\cal P}, \Set]$ in $\cal U$, and for each $X$ in $\cal K$ we have an arrow $\beta_{X}:I(g(X))=([X_{\leq}, \Set], E_{X_{\leq}})\to (\Sh(X), P_{X})=J(X)$ in $\cal U$. It is easy to verify that these arrows satisfy the hypotheses of our method, which yields in this case precisely the Alexandrov coreflection between preorders and topological spaces.
 
\end{enumerate}   

\end{enumerate}

\end{document}